\documentclass[11pt]{article}
\def\WHO{nbd}
\def\version{03.04.2019}\def\users{}  %
\def\users{final-layout}  % when activated, ``our'' debugging is suppressed
\usepackage{amsmath,amsfonts,amssymb,color}
\usepackage{mathrsfs} % THIS PACKAGE ALLOWS FOR ``scr'' FONTS
\usepackage{epsfig}
\usepackage{psfrag}
\usepackage{upgreek} % for upshape greek letters
\usepackage{cite}
\usepackage{url}
\usepackage{import}
\newtheorem{theorem}{Theorem}[section]

\newtheorem{definition}[theorem]{Definition}

\newtheorem{remark}[theorem]{Remark}

\numberwithin{equation}{section}
\usepackage{ifthen}
\ifthenelse{\equal{\users}{final-layout}}{}{
\usepackage{fancyhdr}
\pagestyle{fancy}
\headheight=28pt\headwidth=16cm
\definecolor{gray}{gray}{0.5}
\rhead{\color{gray}visco-elasticity\\
S. Kr\"omer, T.Roub\'\i\v cek}
\chead{}
\lhead{Version \version, file: \jobname.tex, \underline{\Large\color{red}\WHO's working on} \\
compiled:
\number\day.\number\month.\number\year\ at
\the\hour:\ifnum\minute<10 0\fi\the\minute\ h }
}

\definecolor{labelkey}{rgb}{1.,.2,0.}

\newcount\hour \newcount\minute
\hour=\time
\divide \hour by 60
\minute=\time
\loop \ifnum \minute > 59 \advance \minute by -60 \repeat

\usepackage[normalem]{ulem}
\definecolor{brown}{rgb}{0.5,0,0}
\ifthenelse{\equal{\users}{final-layout}}{

    \newcommand{\DELETE}[1]{}
    
    \newcommand{\COMMENT}[1]{}
    
    \newcommand{\COL}[1]{#1}
    \newcommand{\TINY}[1]{}
    \newcommand{\MARGINOTE}[1]{}
}
{

 \newcommand{\DELETE}[1]{{\color{brown}\sout{#1}\color{black}}}
 
 \newcommand{\COMMENT}[1]{{\color{red}\uuline{#1}\color{black}}}
 
 \newcommand{\COL}[1]{{\color{blue}#1}}
 \newcommand{\TINY}[1]{{\tiny#1}}
 \newcommand{\MARGINOTE}[1]{\marginpar{\color{red}\tiny\texttt{#1}}}
}

\renewcommand\dot[1]{\mathchoice
                 {{\buildrel{\hspace*{.1em}\text{\LARGE.}}\over{#1}}}
                 {{\buildrel{\hspace*{.1em}\text{\Large.}}\over{#1}}}
                 {{\buildrel{\hspace*{.1em}\text{\large.}}\over{#1}}}
                 {{\buildrel{\hspace*{.1em}\text{\large.}}\over{#1}}}}
\newcommand\DT{\dot}
\newcommand\DDT[1]{\mathchoice
   {{\buildrel{\hspace*{.1em}\text{\LARGE.\hspace*{-.1em}.}}\over{#1}}}
   {{\buildrel{\hspace*{.1em}\text{\Large.\hspace*{-.1em}.}}\over{#1}}}
   {{\buildrel{\hspace*{.1em}\text{\large.\hspace*{-.1em}.}}\over{#1}}}
   {{\buildrel{\hspace*{.1em}\text{\large.\hspace*{-.1em}.}}\over{#1}}}}

\def\R{{\mathbb R}}

\newcommand\bbD{\mathbb D}

\newcommand\GDir{\mathchoice{\varGamma_{\hspace*{-.15em}\mbox{\tiny\rm D}}}
                          {\varGamma_{\hspace*{-.15em}\mbox{\tiny\rm D}}}
                          {\varGamma_{\hspace*{-.05em}\mbox{\tiny\rm D}}}
                          {\varGamma_{\hspace*{-.05em}\mbox{\tiny\rm D}}}}
\newcommand\GNeu{\mathchoice{\varGamma_{\hspace*{-.15em}\mbox{\tiny\rm N}}}
                          {\varGamma_{\hspace*{-.15em}\mbox{\tiny\rm N}}}
                          {\varGamma_{\hspace*{-.05em}\mbox{\tiny\rm N}}}
                          {\varGamma_{\hspace*{-.05em}\mbox{\tiny\rm N}}}}
\newcommand{\SNeu}{\varSigma_{\mbox{\tiny\rm N}}}
\newcommand{\SDir}{\varSigma_{\mbox{\tiny\rm D}}}
\newcommand{\nablaS}{\nabla_{\scriptscriptstyle\textrm{\hspace*{-.3em}S}}^{}}
\newcommand{\divS}{\mathrm{div}_{\scriptscriptstyle\textrm{\hspace*{-.1em}S}}^{}}
\renewcommand{\d}{{\rm d}}
\newcommand{\DDD}[3]{\begin{array}[t]{c}#1\vspace*{-1em}\\_{#2}\vspace*{-.45em}\\_{#3}\end{array}}
\newcommand{\ddd}[3]{\DDD{\begin{array}[t]{c}\underbrace{#1}\vspace*{.6em}\end{array}}{\text{\footnotesize #2}}{\text{\footnotesize #3}}}
\newcommand{\Vdots}{\mathchoice{\:\begin{minipage}[c]{.1em}\vspace*{-.4em}$^{\vdots}$\end{minipage}\;}{\:\begin{minipage}[c]{.1em}\vspace*{-.4em}$^{\vdots}$\end{minipage}\;}{\:\tiny\vdots\:}{\:\tiny\vdots\:}}
\newcommand{\mathfraks}{\text{\large$\mathfrak{s}$}}
\newcommand{\mathfrakh}{\text{\large$\mathfrak{h}$}}
\newcommand{\scS}{\text{\large$\mathfrak{S}$}}
\newcommand{\barOmega}{{\hspace*{.2em}\overline{\hspace*{-.2em}\varOmega}}}
\newcommand{\barsigma}{{\hspace*{.1em}\overline{\hspace*{-.1em}\sigma}}}
\newcommand{\bary}{{\hspace*{.1em}\overline{\hspace*{-.1em}y}}}
\newcommand{\barf}{{\hspace*{.1em}\overline{\hspace*{-.1em}f}}}

\begin{document}

\MARGINOTE{possible~journals:
\\
\\J. Elast
%\\(or)
%\\MMS,
%\\SIMA,
%\\Nonlin.Anal.,
%\\Nonlinearity (?),
%\\Euro J. Mech. A/Solids (?)
}
\begin{center}
\Huge\bfseries
Quasistatic viscoelasticity with self-contact at large strains
%respecting global non-interpenetration.
\end{center}

\bigskip

\begin{center}
  {\bfseries
  Stefan Kr\"omer}\footnote{\label{fnote1}Institute of
    Information Theory and Automation, Czech Acad. Sci.,
    Pod vod\'arenskou v\v e\v z\'\i\ 4,
     CZ-182 08 Praha 8, Czech Republic}
 and 
{\bfseries Tom\'{a}\v{s} Roub\'\i\v{c}ek}\footnote{Charles University,
Mathematical Institute, Sokolovsk\'a 83, CZ-186~75~Praha~8,
Czech Republic %, {\tt roubicek@karlin.mff.cuni.cz}
}$^,$\footnote{Institute of
Thermomechanics, Czech Acad. Sci., Dolej\v skova~5,
CZ--182~08 Praha 8, Czech Republic}
\end{center}

\bigskip
\noindent{\bf Abstract:} 
The frame-indifferent viscoelasticity in Kelvin-Voigt
rheology at large strains is
formulated in the reference configuration (i.e.\ using the
Lagrangian approach) considering also the possible self-contact
in the actual deformed configuration. Using the concept of
2nd-grade nonsimple materials, existence of certain weak solutions
which are a.e.\ injective is shown by converging an approximate
solution obtained by the implicit time discretisation.

\medskip


\medskip

\noindent{\bf  AMS Clasification:}
%65M12, % PDE evolutionar, Stability and convergence of numerical methods
%65P10, %Numerical problems in dynamical systems, Hamiltonian systems including symplectic integrators
%65Z05, %Applications to physics
%74C10, % Solid mech., Plastic materials, Small-strain
%74F10, % Solid mech., fluid-solid interactions .... porosity,
%74H15, % Solid mech., dynam. probl, Numerical approximation of solutions
%74M15, % Contact
%74R20, % Anelastic fracture and damage
%74S05, % Finite element methods
%76S05. % Flows in porous media; filtration; seepage
 35K86, % Nonlinear parabolic unilateral problems and nonlinear parabolic variational inequalities
 35Q74, % PDEs in connection with mechanics of deformable solids
 74A30, % Nonsimple materials.
 74B20, % Nonlinear elasticity
 74M15. % Contact
 
\section{Introduction}
%        ~~~~~~~~~~~~

Nonlinear elasticity and viscoelasticity is a vital part of
  the continuum mechanics of solids and still faces many open fundamental
  problems even after intensive scrutiny within past many decades. 
  One of such problem is \COL{the possibility of non-physical self-interpenetration
  and analytically supported methods to prevent it}. The problem is difficult
  because of an interaction of two configurations: the reference
  one (ultimately needed for analysis of problems in solid mechanics
  at large strains) and the actual one (ultimately need for determination
  the possible time-varying self-contact boundary region). 

So far, besides merely static situations, only rate-independent evolution
of some internal variables based on (not always very realistic)
concept of instantaneous global minimization and energetic solutions,
possessing a good variational structure and thus allowing incorporation of
the Ciarlet-Ne\v cas condition \cite{CiaNec87ISCN},
has been treated in \cite{MieRou15RIST}. In the viscoelasticity, one cannot
rely purely on a variational structure but should rather work in terms
of partial differential equations.
As emphasized in \cite{MiOrSe14ANVM,MiRoSa18GERV},
``the theory of viscoelasticity at finite strain is notoriously difficult and''
that time it seemed 
``that the present mathematical tools are not sufficient to provide
sufficiently strong solutions in the multidimensional, truly geometrically
invariant case''. Since then, 
the quasistatic viscoelasticity has been treated in \cite{MieRou??TKVR} and
in the dynamical variant in \cite[Sect.\,9.3]{KruRou19MMCM}, but without
globally ruling out self-interpenetration.

In the case of self-contact, instead of differential equations,
it is natural to describe static critical points by variational inequalities.
This was developed for a purely static situation in \cite{PalHea17ISCS} for non-simple materials involving a higher order term in the energy.
For the case of a static obstacle problem neglecting
possible self-contact and self-interpenetration see also \cite{Schu02VACP}.

The goal of this article is to merge the results 
of A.Z.\,Palmer and T.J.\,Healey \cite{PalHea17ISCS} with the
evolution viscoelastic model from \cite{MieRou??TKVR}, using 
a generalization of Korn's inequality developed by P.\,Neff and W.\,Pompe
\cite{Neff02KFIN,Pomp03KFIV}. By this, we obtain  first
results for viscoelastic model allowing a self-contact while
respecting non-self-interpenetration. Let us point out that for a long time, this
was an open problem and largely ignored in engineering
numerical calculations which admitted interpenetration,
relying \COL{solely on the fact that for particular scenarios, 
computational simulations are often not likely} to
go into such non-physical situations.

The plan of the article is following:
In Sect.\,\ref{sect-viscoelast}, we formulate the problem in
terms of the classical partial differential equations, together
with its weak form. Then, in Sect.\,\ref{sect-anal}, 
we employ a time discretisation, prove basic a-priori estimates and,
by convergence for time-step approaching zero, prove existence of
a weak solution. At this point, local non-selfinterpenetration 
and avoidance of singularities while keeping the deformation
gradient ``uniformly'' invertible everywhere is granted by using
the so-called 2nd-grade non-simple (i.e.\ involving strain-gradient)
material concept and the results from \cite{HeaKro09IWSS}.

\section{Quasistatic viscoelasticity}\label{sect-viscoelast}
%        ~~~~~~~~~~~~~~~~~~~~~~~~~~~

\def\pl{\partial}

Strain-gradient theories describe materials referred to as nonsimple, or also multipolar 
or complex. This concept has been introduced a long time ago, cf.\ 
\cite{Toup62EMCS,Toup64TECS} or also e.g.\
\cite{Batr76TNEM,FriGur06TBBC,Podi02CISM,Silh88PTNB,TriAif86GALD}.
%\COMMENT{MinEsh68FSTL was missing = Mindlin, R.D., Eshel, N.N.: On first strain-gradient theories in linear elasticity. Internat. J. Solids Structures 4, 109--124 (1968) - but it is only at small strains, so maybe not urgently to be cited}

We will use the Lagrangian approach and formulate the model in the 
reference (fixed) domain $\varOmega\subset\R^d$
with a smooth boundary $\varGamma:=\pl\varOmega$, $d\ge2$.

To introduce our model in a broader context, we may define  
the {\it total free energy} and the {\it total dissipation potential} 
\begin{subequations}\begin{align}\label{E}
  &\varPsi(y)=\left\{\begin{alignedat}[c]{2}
		&\displaystyle{\int_\varOmega\varphi(\nabla y)+\mathscr{H}(\nabla^2y)\,\d x} \quad &&
		\text{if }\displaystyle{\int_\varOmega\det\nabla y\,\d x
		\leq {\rm meas}\,y(\varOmega),}\\
		&+\infty&&\text{otherwise}
 \end{alignedat}\right.
\intertext{\vspace*{-3ex} and }
	&\mathcal{R}(y,\DT y)=\int_\varOmega\zeta(\nabla y,\nabla\DT y)\,\d x,
\end{align}\end{subequations}
respectively. The condition $\int_\varOmega\det\nabla y\,\d x\le{\rm meas}\,y(\varOmega)$
involved in \eqref{E} is called
the Ciarlet-Ne\v cas condition \cite{CiaNec87ISCN}.
Together with $\det\nabla y\geq 0$ on $\varOmega$ as a property which can
(and will) be ensured by the strain energy $\varphi$, it
guarantees global non-interpenetration.

The mechanical evolution part can then be viewed as an abstract gradient flow 
\begin{align}\label{grad-flow}
\pl_{\DT y}\mathcal{R}(y,\DT y)+\varPsi'(y)
=\mathcal{F}(t)
\ \text{ with }
\big\langle\mathcal{F},y\big\rangle=\int_\varOmega\!
f(x,t){\cdot}y(x) \,\d x, 
\end{align}
cf.\ also \cite{Tved08QEVS,MiOrSe14ANVM} 
for the  case without the Ciarlet-Ne\v{c}as condition. The sum of the
conservative and
the dissipative parts corresponds to the {\it Kelvin-Voigt rheological
  model} in the quasistatic variant (neglecting inertia). Here and henceforth, the notation
``$\,\pl\,$'' is used for partial derivatives (here functional, or
later in Euclidean spaces), while $(\cdot)'$ is used
for the derivative of functions of only one variable.

The generalized gradient $\varPsi'$ is to be understood rather
formally due to the integral Ciarlet-Ne\v cas constraint. This constraint
gives rise to the reaction force due to a possible self-contact.
The contact zone is time evolving and not a-priori known, which is in some
sense a generalization of a so-called Hertz contact at small strains.
At large strains, one must distinguish between the actual deforming
configuration which is relevant for the contact and the reference
configuration which is to be used for analysis and for formulation of
the boundary conditions. Here we use the results of
Palmer and Healey that describe the boundary forces that arise due to the constraint in a static situation \cite{PalHea17ISCS}.

Writing \eqref{grad-flow} locally in the classical formulation, one
arrives at the nonlinear parabolic 4th-order partial differential
equation expressing quasistatic {\it momentum equilibrium},
\begin{align}\label{moment-eq}
\mathrm{div}\,\sigma+g=0 \qquad \text{ with }\quad
\sigma=\sigma_\mathrm{vi}+\sigma_\mathrm{el}
-\mathrm{div}\,\mathfrak{h}_\mathrm{el},
%\quad\text{and}\quad g(y)=g-\pl_y\gamma(y)
\end{align}
where the viscous 
stress is $\sigma_\mathrm{vi}=\sigma_\mathrm{vi}(F,\DT F)$ and the
elastic stress is $\sigma_\mathrm{el}=\sigma_\mathrm{el}(F)$ with $F$ a placeholder
for the deformation gradient $\nabla y$ \COL{and $\dot F$ a placeholder for its time derivative}, while $\mathfrak{h}_\mathrm{el}$ is a
so-called hyperstress arising from the 2nd-grade nonsimple-material concept,
cf.\ e.g.\ \cite{Podi02CISM,Silh88PTNB,Toup62EMCS}.
%% \COMMENT{OK??}   Alex: YES is OK 
In view of the local potentials used in \eqref{grad-flow}, we have
\begin{align}\label{potentials}
&\sigma_\mathrm{vi}(F,\DT F)=\partial_{\DT F}\zeta(F,\DT
  F),
\quad
\sigma_\mathrm{el}(F)=\varphi'(F),
\quad \text{ and} \quad 
\mathfrak{h}_\mathrm{el}(\mathsf{G})=\mathscr{H}'(\mathsf{G}),
\end{align}
\mbox{}where $\mathsf{G}\in \R^{d\times d \times d}$ is a
  placeholder for $\nabla F$, i.e.\ for $\nabla^2y$.
%\COMMENT{WHY $\sfA$?? WHY NOT $\mathsf{G}$ or Hessian $\sfH$??}

\mbox{}%\REPLACE{A physically ultimate}
{An important physical}
requirement is {\it static and dynamic frame indifference}. 
For the elastic stresses, static frame indifference means that 
\begin{subequations}
\begin{align}
  \sigma_\mathrm{el}(RF)= R\, \sigma_\mathrm{el}(F)
  \quad \text{and} \quad  
   \mathfrak{h}_\mathrm{el}(R\mathsf{G})=R\mathfrak{h}_\mathrm{el}(\mathsf{G})
\end{align} 
for all $R\in\mathrm{SO}(d)$, $F$ and $\mathsf{G}$. 
%\COMMENT{Alex: added the conditions that hyperstresses are frame-indifferent!!} 
For the viscous stresses, dynamic frame indifference means that
\begin{align}
\sigma_\mathrm{vi}(RF,\DT RF{+}R\DT
   F)=R\,\sigma_\mathrm{vi}(F,\DT F)
\end{align}
\end{subequations}
for all smoothly time-varying $R:t\mapsto R(t)\in\mathrm{SO}(d)$ \COL{and $F:t\mapsto F(t)\in GL^+(d)$},
cf.\ \cite{Antm98PUVS}.  
%\INSERT
{Note that $R$ may depend on $t$ but
  not on $x\in \varOmega$, since frame-indifference relates to
  superimposing time-dependent \emph{rigid-body motions}.}

In terms of the 
%\INSERT
{thermodynamic} potentials 
%\INSERT{
$\zeta$, $\varphi$, and $\mathscr H$, these frame indifferences read as
\begin{subequations}\label{frame-indif}
\begin{align}
&\varphi(RF)=\varphi(F), \quad 
\mathscr{H}(R G)=\mathscr H(G), 
 \quad \text{and}
\\ 
&\zeta(RF;(RF)\!\DT{^{}}\hspace{.0em})
=\zeta(RF;\DT{R}F{+}R\DT{F})=\zeta(F;\DT{F})
\end{align}
\end{subequations}
for $R$, $F$ and $G$ as above.  
%\DELETE{while no special requirement is needed for $\mathscr{H}$ because the} 
%\DELETE{rotation $R=R(t)$ is independent of $x$ and thus $\nabla(RF)=R\nabla F$.}

As to $\zeta$, the simplest possible choice with such a frame indifference leads to a
viscosity $\sigma_\mathrm{vi}=\partial_{\DT{F}}\zeta$ which is 
linear in $\DT{F}$, while its associated potential is quadratic: 
\begin{align}\label{zeta-special}
  \zeta(F;\DT F):=\frac12 \hat\bbD(C)\DT C{:}\DT C\ \ \ \text{ where }\ \ \ 
  C:=F^\top F \ \text{ and } \ \DT C:=\DT F^\top F+F^\top\DT F.
\end{align}
Notice that frame indifference 
in \eqref{zeta-special} is built in by using a potential which only depends on the right
Cauchy-Green tensor $C$ and its formal time derivative $\DT C$.
To avoid unnecessary technicalities, we adopt this kind of viscosity term for the rest of the paper.
Although the material viscosity is linear as a consequence of \eqref{zeta-special}, the
geometrical nonlinearity arising from large strains is still a vital
part of the problem since $\sigma_\mathrm{vi}(F,\DT F)$ depends on $F$, too.

 Altogether, denoting $\bbD(F):=\hat\bbD(F^\top F)$, we
 arrive at the
 parabolic problem
\begin{align} \label{momentum-eq}
&\mathrm{div}\big(\sigma_\mathrm{vi}(\nabla y,\nabla\DT y)
+\sigma_\mathrm{el}(\nabla
y)-\mathrm{div}\mathscr{H}'(\nabla^2y)\big)+f%(t,y)
=0 
\\ \nonumber
&\qquad\qquad\qquad\qquad\qquad \text{with }\ \
\sigma_\mathrm{vi}(F,\DT F){:}= 2F\bbD(F)(F^\top\DT F{+}\DT F^\top F)
\\ \nonumber
&\qquad\qquad\qquad\qquad\qquad \text{and }\ \
\sigma_\mathrm{el}(F)=\varphi'(F)\,,
\end{align}
on $Q$. 

We complete \eqref{momentum-eq} by some boundary conditions. 
For simplicity, 
we only consider
a mechanically fixed part $\GDir$, undeformed and independent of time (i.e.\ identity): 
\begin{subequations}
  \label{BC}
 \begin{align}
&y(x)=x\ \ \ \text{(identity)}&&\text{on }\ \GDir, \label{BC2}
\\&\big(\sigma_\mathrm{vi}(\nabla y,\nabla\DT y)
+\sigma_\mathrm{el}(\nabla y)\big)\vec{n}
-\divS(\mathscr{H}'(\nabla^2y)\vec{n})=\mathfraks\qquad\text{(a reaction\ force)}&&\text{on }\ \GNeu, \label{BC1} 
\\&\label{BC3} 
\mathscr{H}'(\nabla^2y){:}(\vec{n}\otimes\vec{n})=0&&\text{on }\ \varGamma,
\end{align}
\end{subequations}
where $\varGamma=\partial\varOmega$, $\GNeu=\varGamma \setminus \GDir$ and $\vec n$ is the outward pointing normal vector.
Moreover,
``$\divS$'' in \eqref{BC1} denotes the surface divergence defined as
$\divS(\cdot)=\mathrm{tr}\big(\nablaS(\cdot)\big)$, 
{where $\mathrm{tr}(\cdot)$ denotes the trace and $\nablaS$
denotes the surface gradient given by $\nablaS v=(\mathbb I- \vec
n{\otimes}\vec n)\nabla v= \nabla v-\frac{\partial
  v}{\partial\vec{n}}\vec{n}$.}
Note that for equilibria, \eqref{BC1} and \eqref{BC3} reduce to the natural boundary conditions complementing the Dirichlet condition  
\eqref{BC2}.

The energetics of the system \eqref{momentum-eq}--\eqref{BC} can be revealed 
by testing \eqref{momentum-eq} by $\DT y$, and using 
the boundary conditions after integration over $\varOmega$ and 
using Green's formula twice together with another
$(d{-}1)$-dimensional Green formula over $\varGamma$ for \eqref{momentum-eq}.
The last mentioned technique is related with the concept of
nonsimple materials; for the details about how the boundary conditions
are handled see e.g.\ \cite[Sect.\,2.4.4]{Roub13NPDE}.
This test of \eqref{momentum-eq} gives the mechanical energy balance:
\begin{align}\label{mech-engr}
\int_\varOmega\!\!\!\! \ddd{2\zeta(\nabla y,\nabla\DT y)_{_{_{_{}}}}}{dissipation}{rate}\!\!\!\!
+\!\!\!\!\ddd{\sigma_\mathrm{el}(\nabla y){:}\nabla\DT y_{_{_{_{}}}}}{mechanical}{power}\!\!\!\!\d x
+\frac{\d}{\d t}
\int_\varOmega\!\!\!\!\!\!\!\!\!\!\ddd{\mathscr{H}(\nabla^2y)_{_{_{_{}}}}}{``nonsimple'' part of}{the stored energy}\!\!\!\!\!\!\!\!\!\!\,\d x=\int_\varOmega\!\!\!\!\!\!\!\!\!\!\ddd{f\cdot\DT y_{_{_{_{}}}}}{power of the}{bulk force}
%{mechanical load}
\!\!\!\!\!\!\!\!\!\!\d x\,.
%+\int_{\GNeu}\!\!\!\!\!\!\ddd{.........\cdot\DT y_{_{_{_{}}}}}{power of}{the traction}\!\!\!\!\!\d S.
\end{align}

In what follows, we will use the (standard) notation for the Lebesgue 
$L^p$-spaces and $W^{k,p}$ for Sobolev spaces whose $k$-th distributional 
derivatives are in $L^p$-spaces and the abbreviation $H^k=W^{k,2}$. 
The notation $W^{1,p}_{\mathrm D}$ will indicate the closed subspace
of $W^{1,p}$ with zero traces on $\GDir$. The Banach space of continuous
functions on a compact set will be denoted as $C(\cdot)$, while their dual
as ${\rm Meas}(\cdot)$ being the space of finite Radon measures; if
scalar-valued, the subset of non-negative measures will be denoted
by ${\rm Meas}^+(\cdot)$.
Moreover, we will use the standard notation  $p'=p/(p{-}1)$.
%, and $p^*$ for the Sobolev exponent $p^*=pd/(d{-}p)$ for $p<d$ while
%$p^*<\infty$ for $p=d$ and $p^*=\infty$ for $p>d$.
%, and the ``trace exponent'' $p^\sharp$ defined 
%as $p^\sharp=(pd{-}p)/(d{-}p)$ for $p<d$ while
%$p^\sharp<\infty$ for $p=d$ and $p^\sharp=\infty$ for $p>d$.
%\COMMENT{Alex 20.7.17: If $p^*$ is used, then it should be defined!!}
%Thus, e.g., $W^{1,p}(\varOmega)\subset L^{p^*}\!(\varOmega)$ or 
%$L^{{p^*}'}\!(\varOmega)\subset W^{1,p}(\varOmega)^*$=\,the dual to $W^{1,p}(\varOmega)$. 
In the vectorial case, we will write $L^p(\varOmega;\R^n)\cong L^p(\varOmega)^n$ 
and $W^{1,p}(\varOmega;\R^n)\cong W^{1,p}(\varOmega)^n$.
%
%Thus, for example,
%NOT USED: Moreover, we denote
%\begin{align}
%H^{1}_{\mathrm D}(\varOmega;\R^{d})=\big\{v\in L^2(\varOmega;\R^d);\ \nabla v\in L^2(\varOmega;\R^{d\times d}),\ \ v|_{\GDir}^{}=0\,\big\}.
%\end{align}
For the fixed time interval $I=[0,T]$, we denote by $L^p(I;X)$ the 
standard Bochner space of Bochner-measurable mappings $I\to X$ with $X$ a 
Banach space whose norm is in $L^p(I)$. Also, $W^{k,p}(I;X)$ denotes the
Banach space of mappings 
from $L^p(I;X)$ whose $k$-th distributional derivative in time is also in 
$L^p(I;X)$. The dual space to $X$ will be denoted by $X^*$.
%Moreover, we denote by $\mathrm{BV}(I;X)$ the Banach space
%of the mappings $I\to X$ that have bounded variation on $I$.
%By $\mathrm{Meas}(I;X)$ we denote the space of $X$-valued measures on $I$. 
Moreover, $C_{\rm w}(I;X)$ denotes the Banach space of weakly continuous
functions $I\to X$, and $L^\infty_{\rm w}(I;X)$ denotes the Banach space of 
essentially bounded, weakly measurable functions $I\to X$. 
The scalar product between vectors, matrices, or 3rd-order tensors
will be denoted by ``$\,\cdot\,$'', ``$\,:\,$'', or ``$\,\Vdots\,$'', 
respectively. Finally, in what follows, $K$ denotes a positive, possibly 
large constant.

We consider an initial-value problem, imposing the initial conditions
\begin{align}\label{IC}
y(0,\cdot)=y_0\ \ \ \ \text{ on }\ \varOmega.
\end{align}

\begin{definition}[Weak solution]\label{def}
 The couple $(y,\mathfraks)$ of a displacement field $y:Q\to\R^d$
 and a reaction traction $\mathfraks$ as a
 distribution $\SNeu\to\R^d$  is called a weak solution
 of the constrained initial-boundary-value problem
\eqref{momentum-eq}--\eqref{BC}--\eqref{IC} 
%\DELETE{, if the initial conditions \eqref{IC} hold,} 
if the following three conditions are satisfied:
\begin{itemize}
\item[(i)]
$(y,\mathfraks)\in %L^\infty 
C_{\rm w}(I;W^{2,p}(\varOmega;\R^d))\times
L^2(I;W^{2-1/p,p}(\GNeu;\R^d)^*)$
%{\rm meas}^+(\SNeu)$
with $\nabla\DT y\in L^2(Q;\R^{d\times d})$ 
%\TT with 
%${\rm div}^2(\mathscr{H}'(\nabla y^2))\in L^2(I;H^1(\varOmega;\R^d)^*)$, \EE
and with $\min_Q\det\nabla y>0$ and $y|_{\SDir}=\text{\rm identity}$, and
%it satisfies
the integral identity 
\begin{subequations}\begin{align}\nonumber
&\int_Q \bbD(\nabla y)(\nabla\DT y^\top\nabla y{+}\nabla y^\top\nabla\DT y){:}(\nabla
y^\top\nabla z+\nabla z^\top\nabla y)+\varphi'(\nabla y)
{:}\nabla z
\\[-.5em]&\qquad\qquad\qquad\qquad\qquad
+\mathscr{H}'(\nabla^2y)\Vdots\nabla^2z\,
\d x\d t=\int_Q f{\cdot}z\,\d x\d t
+
\big\langle\mathfraks,z|_{\SNeu}^{}\big\rangle
\label{momentum-weak-def}
\end{align}
is satisfied
for all 
smooth $z:Q\to\R^d$ with $z=0$ on $\SDir$ together with $y(0,\cdot)=y_0$. 
\item[(ii)] For a.e.~$t\in I$, $y(t,\cdot)$ satisfies the Ciarlet-Ne\v{c}as condition, i.e.,
\begin{align} \label{CiarletNecas}
  \int_\varOmega\det\nabla y(t,x)\,\d x\le{\rm meas}\,y(t,\varOmega).
\end{align}
\item[(iii)] The support of $\mathfraks$ is contained in the part of $\SNeu$
which, after deformation, is in self-contact: 
\begin{align} \label{PalHea-traction-limit}
  \big\langle\mathfraks(t,\cdot),z\big\rangle=0 \quad
  \text{for a.e.~$t\in I$ and every $z\in W^{2-1/p,p}(\GNeu;\R^d)$ vanishing on $\scS_t$},
\end{align}\end{subequations}
where the self-contact set at time $t$ given by
\[
	\scS_t:= \{ x\in \GNeu \mid \exists \tilde{x}\in \barOmega\setminus \{x\}:\ y(t,x)=y(t,\tilde{x})\}.
\]
\end{itemize}
\end{definition}

\begin{remark}[The role of $\mathfraks$]\upshape
  The constraint \eqref{CiarletNecas}, effectively ruling out
  self-interpenetration, is also built into our definition \eqref{E}
  of the total free energy $\Psi$. The reaction force $\mathfraks$ has the role of a Lagrange multiplier which only appears
  at the ``boundary'' of this constraint, i.e., when there is self-contact.
\end{remark}

\begin{remark}[Frame indifference and more general viscosity terms]\upshape
Even in the case when $\zeta$ is not quadratic, the frame indifferences \eqref{frame-indif} imply (cf.~\cite{Antm95NPE,MiOrSe14ANVM})
the existence of 
reduced potentials $\hat\varphi$,  $\hat\zeta$, and $\hat{\mathscr H}$  such that
\begin{align}\label{hat-ansatz}
\zeta(F,\DT F)=\hat\zeta(C,\DT C), \quad 
\varphi(F)=\hat\varphi(C),  \quad \text{and}\quad
\mathscr H(\mathsf{G}) = \hat{\mathscr H}(\mathsf B) 
\end{align}
where $\mathsf B=\mathsf{G}^\top\!\cdot \mathsf{G} \in \R^{(d\times d)\times (d
  \times d)} $, and $C\in \R^{d \times d}_\mathrm{sym}$ is the right
Cauchy-Green tensor 
with time derivative $\DT C$ as in \eqref{zeta-special}.
More specifically, denoting $\mathsf{G}=[\mathsf{G}_{\alpha ij}]$ the placeholder 
for $\frac{\pl}{\pl x_j}F_{\alpha i}$ with $F_{\alpha i}$ the placeholder 
for $\frac{\pl}{\pl x_i}y_{\alpha}$, the exact meaning is 
$[\mathsf{G}^\top\!\cdot \mathsf{G}]_{ijkl}:=\sum_{\alpha=1}^d\mathsf{G}_{\alpha ij}\mathsf{G}_{\alpha kl}$
and $[F^\top F]_{ij}:=\sum_{\alpha=1}^dF_{\alpha i}F_{\alpha j}$.
%\REPLACE{This}
{The ansatz \eqref{hat-ansatz}} also means that
\begin{subequations}
\label{KV-large-F-vs-C}
\begin{align}
&\sigma_\mathrm{el}=2 F\pl_C^{}\hat\varphi(F^\top F), \quad 
\mathfrak{h}_\mathrm{el}(\mathsf{G}) = 2\mathsf{G} \pl_{\mathsf B} \hat{\mathscr
  H}(\mathsf{G}^\top\!\!\cdot \mathsf{G}) = 2\mathsf{G} \pl_{\mathsf B} \hat{\mathscr
  H}(\mathsf B), \\
&\sigma_\mathrm{vi}(F,\DT F)
=2 F\pl_{\DT C}\hat\zeta(F^\top F,\DT F^\top F{+}F^\top\DT F)
=2 F\pl_{\DT C}\hat\zeta(C,\DT C).
\end{align}
\end{subequations}
Furthermore, the \emph{specific dissipation rate} can easily be identified 
in terms of $\hat\zeta$ as 
\begin{align}\nonumber
\xi(F,\DT F)&=\sigma_\mathrm{vi}(F,\DT F){:}\DT F=
2F\pl_{\DT C}\hat\zeta(F^\top F,\DT F^\top F{+}F^\top \DT F){:}\DT F
\\&=
 \pl_{\DT C}\hat\zeta(F^\top F,\DT F^\top F{+}F^\top \DT F){:}
(\DT F^\top F{+}F^\top \DT F).
\label{def-of-xi}\end{align}
For our choice \eqref{zeta-special}, we thus have 
$\xi(F,\DT F)=\hat\bbD(C)\DT C{:}\DT C
=2\zeta(F;\DT F)$. 
\end{remark}

\section{Analysis by time discretisation}\label{sect-anal}
%        ~~~~~~~~~~~~~~~~~~~~~~~~~~~~~~~

Let us summarize the assumptions we impose on the data:
\begin{subequations}\label{ass}
\begin{align}
&\nonumber\exists\, p>d,\; \COL{s>1},\; q\ge pd/(p{-}d)\ \exists\,
\alpha,K,\epsilon>0: 
\\
\nonumber&\varphi:\mathrm{GL}^+(d)\to\R^+\ \text{ continuously 
differentiable},\ \forall\, F\in\mathrm{GL}^+(d):
\\&\label{ass-phi}
\qquad
\varphi(F)\ge\epsilon|F|^s+\epsilon/(\det F)^q,
\\[0.4em]
&\nonumber\mathscr{H}:\R^{d\times d\times d}\to\R^+
\text{convex, continuously differentiable}, \ \forall\, G,G_1,G_2 \in\R^{d\times d\times d}:
\\
&\label{ass-H1}
\qquad
\epsilon|G|^p\le\mathscr{H}(G) \leq K(1{+}|G|^p),\ |\mathscr{H}'(G)|\leq K(1{+}|G|^{p-1}),
\\[0.4em]
&\label{ass-H2}
\qquad
\alpha |G_1-G_2|^p\leq [\mathscr{H}'(G_1)-\mathscr{H}'(G_2)]\Vdots(G_1-G_2),
\\[0.4em]
&\nonumber
\zeta(F,\DT F)=\textstyle{\frac12} \bbD(F) (F^\top \DT F+ \DT F^\top F):(F^\top \DT F+ \DT F^\top F)
\ \text{with $\bbD(F)=\hat \bbD(F^\top F)$ (cf.\ \eqref{zeta-special})}, \\
&\nonumber 
\text{$C\mapsto \hat\bbD(C)$ continuous and bounded, } 
\forall\, C, \DT C \in \R^{d \times d}_\mathrm{sym}: \\
&\qquad  \label{ass-D1}
\hat\bbD(C):\R^{d \times d}_\mathrm{sym} \to \R^{d \times d}_\mathrm{sym} \ \ \text{is linear and symmetric,} \\
&\qquad  \label{ass-D2}
\hat\bbD(C)\DT C:\DT C \geq \alpha |\DT C|^2,
\\[0.4em]
&\label{ass-f}
f\!\in\! H^1(I;L^2(\varOmega;\R^d)),
\\[0.3em]
&\label{ass-IC}
y_0\in W^{2,p}(\varOmega;\R^d), \quad \mathrm{det}(\nabla y_0)\ge\epsilon,\quad y_0|_{\GDir}=\text{identity},
\end{align}\end{subequations}
where $\mathrm{GL}^+(d)$ denotes the set of matrices in $\R^{d\times  d}$ with
positive determinant.

The balance of the mechanical energy \eqref{mech-engr} can be written in
the more specific form 
\begin{align}%\nonumber
  &\int_\varOmega\!\!\!\!
  \ddd{\zeta(\nabla y,\nabla\DT y)_{_{_{_{}}}}}{dissipation}{rate}\!\!\!\!\!\d x
+\frac{\d}{\d t}
\int_\varOmega\!\!\!\!\!
\ddd{\varphi(\nabla y)+\mathscr{H}(\nabla^2y)_{_{_{_{}}}}}{stored}{ energy}\!\!\!\!\!\d x
%\\& \qquad\qquad 
=\int_\varOmega\!\!\!\!\!
\ddd{f{\cdot}\DT y_{_{_{_{}}}}}{power of}{bulk load}\!\!\!\!\!\!\!\!
\d x\,.
\label{mech-engr+}\end{align}
Our main result is the following:

\begin{theorem}[Existence of  weak solutions]
\label{thm:MainExist}
Let \eqref{ass} hold. 
Then there exists a weak solution to the constrained initial-boundary-value problem
\eqref{momentum-eq}--\eqref{BC}--\eqref{IC} in the sense of Definition~\ref{def}.
\end{theorem}

%\begin{proof}
\noindent{\it Proof.}
As we have neglected inertial effects, we can use time discretisation.
We consider a time step $\tau>0$ such that $T/\tau$ is integer, 
and an equidistant partition of the time interval $I=[0,T]$.
Thus the regularized system \eqref{momentum-eq}--\eqref{BC}
after this discretisation takes the form 
\begin{subequations}\label{system-disc}
\begin{align}\nonumber
&-\mathrm{div}\bigg(\sigma_\mathrm{vi}\Big(\nabla y_\tau^{k-1},
\nabla\frac{y_\tau^k{-}y_\tau^{k-1}\!\!}\tau\ \Big)
+\sigma_\mathrm{el}(\nabla y_\tau^k)
\\[-.3em]&\qquad\qquad\qquad\qquad\qquad
-\mathrm{div}\mathscr{H}'(\nabla^2y_\tau^k)\ \bigg)
=f_\tau^k:=\frac1\tau\int_{(k-1)\tau}^{k\tau}\!\!f(t)\,\d t&&\text{on }\ \varOmega,
\\
%\nonumber
&\bigg(\sigma_\mathrm{vi}\Big(\nabla y_\tau^{k-1},\nabla\frac{y_\tau^k{-}y_\tau^{k-1}\!\!}\tau\ \Big)
+\sigma_\mathrm{el}(\nabla y_\tau^k)
\bigg)\vec{n}%\\&\qquad\qquad\qquad\qquad\qquad\qquad
-\divS(\mathscr{H}'(\nabla^2y_\tau^k)\vec{n})=\mathfraks_\tau^k&&\text{on }\ \GNeu, 
\\&y_\tau^k=\text{\,identity\ \ on }\ \GDir,\qquad\qquad \mathscr{H}'(\nabla^2y_\tau^k){:}(\vec{n}\otimes\vec{n})=0&&\text{on }\ \varGamma,
\label{BC1-reg-disc}
\end{align}\end{subequations}
which is to be solved recursively for $k=1,...,T/\tau$, starting with
$y_\tau^0=y_0$. 

This boundary-value problem (in its suitable weak formulation) has a
variational structure. More specifically, a weak solution can be obtained
from the problem:
\begin{align}\label{E-R}
  &\left.\begin{array}{ll}\text{minimize}&
\displaystyle{\int_\varOmega \varphi(\nabla y)+ \mathscr{H}(\nabla^2y)
  +\tau\zeta\Big(\nabla y_\tau^{k-1}, \frac{\nabla y{-} \nabla y_\tau^{k-1}}{\tau}   %\nabla\frac{y{-}y_\tau^{k-1}}\tau
  \Big)-f_\tau^k\cdot y\,\d x}
\\[.0em]\text{subject to}&
\displaystyle{{\rm meas}\,y(\varOmega)\ge\int_\varOmega\det\nabla y\,\d x\,,\ \ 
y=0~\text{on}~\GDir, \ \
  y\in W^{2,p}(\varOmega;\R^d)\,.}
\end{array}\right\}
\end{align}
By the standard direct-method arguments, this problem has a solution which
we will denote as $y_\tau^k$.

Comparing the value of the functional in the first line of \eqref{E-R} at $y=y_\tau^k$ with its value
at $y=y_\tau^{k-1}$ which must be bigger of equal, and summing it for
$k=1,...,K$, we obtain the discrete energy imbalance for arbitrary $K\leq T/\tau$:
\begin{align}
\begin{aligned}\label{eimba}
  &\int_\varOmega\varphi(\nabla y_\tau^K)+\mathscr{H}(\nabla^2y_\tau^K)\,\d x
  +\tau\sum_{k=1}^K \int_\varOmega\zeta\Big(\nabla y_\tau^{k-1},%\nabla y{-}\nabla y_\tau^{k-1}
  \nabla\frac{y_\tau^k{-}y_\tau^{k-1}}\tau\Big)\,\d x
\\[-.6em]%\nonumber
&\quad\le
\tau\sum_{k=1}^K \int_\varOmega f_\tau^k\cdot\frac{y_\tau^k{-}y_\tau^{k-1}}\tau\,\d x
+\int_\varOmega\varphi(\nabla y_0)+\mathscr{H}(\nabla^2y_0)\,\d x.
\\[-.6em]%\nonumber
&\quad=\int_\varOmega f_\tau^K\cdot y_\tau^K \,\d x
-\tau\sum_{k=1}^K \int_\varOmega y_\tau^{k-1}\cdot\frac{f_\tau^k{-}f_\tau^{k-1}}\tau\,\d x
+\int_\varOmega\varphi(\nabla y_0)+\mathscr{H}(\nabla^2y_0)
-f_\tau^0\cdot y_0 \,\d x\,
\\[-.6em]%\nonumber
&\quad \COL{
\leq C\|f\|_{H^1(I;L^2(\varOmega;\R^n))} \sup_{0\leq k\leq T/\tau} \|y_\tau^{k}\|_{L^2(\varOmega;\R^n)}
+\int_\varOmega\varphi(\nabla y_0)+\mathscr{H}(\nabla^2y_0)
-f_\tau^0\cdot y_0 \,\d x\,.
}
\end{aligned}
\end{align}

Considering $\{y_\tau^k\}_{k=0,...,T/\tau}$,
we introduce a notation for the piecewise-constant and the piecewise affine 
interpolants defined respectively by
\begin{subequations}\label{def-of-interpolants}
\begin{align}\label{def-of-interpolants-}
&&&
\bary_\tau(t)= y_\tau^k\,,\qquad\ \
\underline y_\tau(t)= y_\tau^{k-1},\qquad\ \
%\underline{\bary}_\tau(t)=y_\tau^{k-1/2},
&&\text{and}
\\&&&\label{def-of-interpolants+}
y_\tau(t)=\frac{t-(k{-}1)\tau\!\!}\tau\  y_\tau^k+\frac{k\tau-t}\tau y_\tau^{k-1}
&&\text{for }(k{-}1)\tau<t\le k\tau.&&&&
\end{align}\end{subequations}
We will also use the notation $\barsigma_\tau$ and $\barf_\tau$ with analogous meaning.

\COL{
Since $\zeta\geq 0$, taking the supremum over $K$ in \eqref{eimba} and
using the Poincar\'e inequality based on the Dirichlet condition and the coercivity \eqref{ass-H1} and \eqref{ass-phi} 
of $\mathscr{H}$ and $\varphi$, respectively, 
%\eqref{ass-D2}: $\zeta$, and the discrete Gronwall inequality, 
we obtain the a-priori estimate
\begin{subequations}\label{est}\begin{align}\label{est-y}
 &\|y_\tau\|_{L^\infty(I;W^{2,p}(\varOmega;\R^d))}\le C\,,
%\\&\|(\nabla \overline{y}_\tau)^\top \nabla\DT y_\tau + (\nabla\DT y_\tau)^\top \nabla \overline{y}_\tau \|_{L^2(Q;\R^{d\times d})}\le C\,\label{est-sym-grad-y}
%\\&\|\nabla\DT y_\tau\|_{L^2(Q;\R^{d\times d})}\le C\,;\label{est-grad-y}
 \intertext{with some constant $C=C(I,f,y_0,\varOmega,\GDir,p,s,d)>0$. Using \cite{HeaKro09IWSS},
   from \eqref{est-y} and \eqref{ass-phi} we can also deduce that $\det\nabla y_\tau>0$ and even that}
    &\Big\|\frac1{\det\nabla y_\tau}\Big\|_{L^\infty(Q)}\le C. \label{est-det}
 \intertext{In addition, by a variant of Korn's inequality \cite{Pomp03KFIV} (cf.~\eqref{PKorn} below) and \eqref{est-y},  
we can also exploit the coercivity \eqref{ass-D2} of $\zeta$ in \eqref{eimba}, which gives that}
	&\|\nabla\DT y_\tau\|_{L^2(Q;\R^{d\times d})}\le C\,.\label{est-grad-y}
\end{align}\end{subequations}
More precisely, for the proof of \eqref{est-grad-y} we used the following generalized Korn inequality proved
by W.\,Pompe \cite{Pomp03KFIV}, generalizing earlier results by P.\,Neff \cite{Neff02KFIN}:
\begin{align} \label{PKorn}
	\|\nabla\DT y_\tau\|_{L^2(Q;\R^{d\times d})}
	\le C \|F \nabla \DT y_\tau + (\nabla\DT y_\tau)^\top F^\top \|_{L^2(Q;\R^{d\times d})},
\end{align}
for a field 
$F\in C(\barOmega;\R^{d\times d})$ with $\min_{\barOmega}^{}\det\,F>0$, 
here applied in a further generalized form with $F:=(\nabla \bary_\tau)^\top$, which is always contained in a fixed compact subset of the admissible fields $F$ due to the uniform bounds \eqref{est-y} and \eqref{est-det}.
}

By the results from \cite{PalHea17ISCS} applied to \eqref{E-R},
we can claim that $y_\tau^k\in W^{2,p}(\varOmega;\R^d)$
satisfies also the identity
\begin{align}\nonumber
&\int_\varOmega \bbD(\nabla y_\tau^{k-1})([\nabla\DT
y_\tau^k]^\top\nabla y_\tau^{k-1}{+}[\nabla y_\tau^{k-1}]^\top\nabla\DT y_\tau^k){:}([\nabla
y_\tau^{k-1}]^\top\nabla z
+[\nabla z]^\top \nabla y_\tau^{k-1}
)+\varphi'(\nabla y_\tau^k)
%\big)
{:}\nabla z%\,\d x
\\[-.5em]&\qquad\qquad\qquad\qquad\qquad
+\mathscr{H}'(\nabla^2y_\tau^k)\Vdots\nabla^2z\,\d x=\int_\varOmega f_\tau^k{\cdot}z\,\d x
+\int_{\GNeu}(\nabla y_\tau^k)^{-\top}\vec{n}{\cdot}z\,\d \sigma_\tau^k
%\big\langle\mathfraks,z\big\rangle
%?????????????
\label{momentum-weak}
\end{align}
for all $z\in W^{2,p}(\varOmega;\R^d)$ with $z|_{\GDir}^{}=0$ and with some scalar-valued non-negative
measure $\sigma_\tau^k\in{\rm Meas}^+(\GNeu)$.
Here, note that to apply \cite{PalHea17ISCS}, we have temporarily
interpreted the $\zeta$-term as absorbed into the elastic energy density $\varphi$ with $\tau>0$ and $\nabla y_\tau^{k-1}\in C(\barOmega;\R^{d\times d})$ fixed; this
possibly breaks frame indifference which is assumed but not exploited in \cite{PalHea17ISCS}.
The expression $\mathfraks_\tau^k:=(\nabla y_\tau^k|_{\GNeu}^{})^{-\top}\vec{n}\sigma_\tau^k$ occurring in the last integral
was obtained in \cite{PalHea17ISCS}. %Note that $(\nabla y_\tau^k|_{\GNeu}^{}){-\top}\vec{n}$ is a vector in direction of the outer normal in the actual deformed configuration.
It is in the position of a traction \COL{in direction of the outer normal in the actual deformed configuration which (up to a positive scalar factor)
is given by $(\nabla y_\tau^k|_{\GNeu}^{})^{-\top}\vec{n}$. Also notice that both $\nabla y_\tau^k$ and its inverse $(\nabla y_\tau^k)^{{-}1}$ (by \eqref{est-det})
are uniformly bounded and (even H\"older) continuous on the closure of $\varOmega$, and so are their traces on $\GNeu$. 
In particular, the traction $\sigma_\tau^k$ itself is a measure. }

As shown in \cite{PalHea17ISCS}, $\sigma_\tau^k$ therefore $\mathfraks_\tau^k$ vanishes outside the self-contact set, i.e.,
\begin{align}
\label{PalHea-traction}
0=\sigma_\tau^k(\{x\in \varGamma_N \mid y_\tau^k(t,x)\neq y_\tau^k(\tilde{x}) \ \text{for all} \ \tilde{x}\in \varGamma_N \setminus\{x\} \}).
\end{align}
%In addition, it also satisfies an action-reaction principle (the traction forces on the two different sides in a contact point sum up to zero).

By comparison, we obtain an estimate on the measure $(\nabla y_\tau^k)\vec{n} \sigma_\tau^k\in {\rm Meas}(\GNeu;\R^d)$,
but unfortunately in a bigger space than the space of measures. Namely, writing \eqref{momentum-weak} in terms of the interpolants as
\begin{align}\nonumber
&\int_Q \bbD(\nabla\underline y_\tau)([\nabla\DT
    y_\tau]^\top\nabla\underline{y}_\tau{+}[\nabla\underline{y}_\tau]^\top\nabla\DT y_\tau){:}([\nabla\underline{y}_\tau]^\top\nabla z+[\nabla z]^\top\nabla \underline{y}_\tau)
  +\varphi'(\nabla\bary_\tau){:}\nabla z%\,\d x\d t
\\[-.5em]&\qquad\qquad\qquad\qquad\qquad
+\mathscr{H}'(\nabla^2\bary_\tau)\Vdots\nabla^2z\,\d x\d t=\int_Q \barf_\tau{\cdot}z\,\d x\d t
+\int_{\SNeu}(\nabla\bary_\tau)^{{-}\top}\vec{n}{\cdot}z\,\d \barsigma_\tau
%\big\langle\mathfraks,z\big\rangle
%?????????????
\label{momentum-weak+}
\end{align}
for all $z\in L^1(I;W^{2,p}(\varOmega;\R^d))$ with $z|_{\SDir}^{}=0$,
we can estimate
\begin{align}\nonumber
&\sup_{\|z\|_{L^2(I;W^{2-1/p,p}(\GNeu;\R^d))}\le1}\int_{\SNeu}(\nabla\bary_\tau)^{{-}\top} \vec{n}{\cdot}z\,\d \barsigma_\tau
\\&\nonumber\qquad\qquad
\leq C_{1} \sup_{\|z\|_{L^2(I;W^{2,p}(\varOmega;\R^d))}\le1}\int_{\SNeu}(\nabla\bary_\tau)^{{-}\top}\vec{n}{\cdot}z|_{\SNeu}^{}\,\d \barsigma_\tau
\\&\nonumber\qquad\qquad
=C_{1} \sup_{\|z\|_{L^2(I;W^{2,p}(\varOmega;\R^d))}\le1}\int_Q \bbD(\nabla\underline y_\tau)([\nabla\DT
    y_\tau]^\top\nabla\underline{y}_\tau{+}[\nabla\underline{y}_\tau]^\top\nabla\DT y_\tau){:}([\nabla\underline{y}_\tau]^\top\nabla z+[\nabla z]^\top\nabla \underline{y}_\tau)
\\[-.7em]&\nonumber\qquad\qquad\qquad\qquad\qquad\qquad\qquad\qquad\qquad\qquad
+\varphi'(\nabla\bary_\tau){:}\nabla z
+\mathscr{H}'(\nabla^2\bary_\tau)\Vdots\nabla^2z-\barf_\tau{\cdot}z\,\d x\d t
\\[.2em]&\nonumber\qquad\qquad
\le %\begin{aligned}[t]
C_{2}\Big((\max|\bbD|) \big\|\nabla\underline y_\tau\|_{L^\infty(Q;\R^{d\times d})} 
\big\|\nabla\DT y_\tau\|_{L^2(Q;\R^{d\times d})}^{}+\big\|\varphi'(\nabla\bary_\tau)\big\|_{L^\infty(Q;\R^{d\times d})}^{}
\\&\qquad\qquad\qquad\qquad\qquad\qquad\qquad
+\big\|\mathscr{H}'(\nabla^2\bary_\tau)\|_{L^2(I;L^{p'}(\varOmega;\R^{d\times d\times d}))}^{}
+\|f\|_{L^2(I;L^2(\varOmega;\R^{d}))}\Big)
%\end{aligned}
\label{est-traction-}
\end{align}
with some constants $C_{1}$, $C_{2}$ depending on $\varOmega$, $d$ and $p$. Together with \eqref{est},
this implies the estimate
\begin{align}\label{est-traction}
\big\|(\nabla\bary_\tau)^{{-}\top}|_{\SNeu}^{} \vec{n}\barsigma_\tau\big\|_{L^2(I;W^{2-1/p,p}(\GNeu;\R^d)^*)}^{}\le C\,.
\end{align}

By the Poincar\'e inequality, \eqref{est-grad-y} together with the time-constant Dirichlet
boundary conditions gives even the estimate on $\DT y_\tau$ in $L^2(I;H^1(\varOmega;\R^d))$.

The estimates (\ref{est}a,c) and \eqref{est-traction} hold for the piecewise constant
interpolants $\bary_\tau$ and $\underline y_\tau$, as well.
Therefore, now we can select a subsequence converging for $\tau\to 0$ in the sense
\begin{subequations}\label{conv}
	\begin{align}\label{conv1}
          &y_\tau\to y,\ \
          %\text{ and }\ \
          \bary_\tau\to y,\ \ \text{ and }\ \ \underline y_\tau\to y
    &&\text{weakly* in }L_{\rm w}^\infty(I;W^{2,p}(\varOmega;\R^d))\,,
    \\&\DT y_\tau\to \DT y&&\text{weakly in }\ L^2(I;H^1(\varOmega;\R^d))\,,
    \\& (\nabla\bary_\tau|_{\GNeu}^{})\vec{n}\barsigma_\tau\to\mathfraks_1 && \text{weakly in }\ L^2(I;W^{2-1/p,p}(\GNeu;\R^d)^*)\,.  \label{conv3}
    \intertext{Note that the limit $y$ also inherits \eqref{est-det} and	
the Ciarlet-Ne\v{c}as condition \eqref{CiarletNecas} from $\bary_\tau$,
since for a.e.~$t\in I$, $\bary_\tau\to y$ in $C^1$ and
${\rm meas}\,\bary_\tau(\varOmega)\to{\rm meas}\,y(\varOmega)$ (see \cite[Prop.~4.3]{MieRou16a}, e.g.).			
		Moreover, by the Aubin-Lions compact-embedding theorem (see \cite[Lemma 7.7]{Roub13NPDE}) and its generalization for time derivative
      measures (see \cite[Cor.\,7.9]{Roub13NPDE}), respectively, we also have that}
    &y_\tau\to y,\ \ \ \bary_\tau\to y,\ \ \text{ and }\ \ \underline y_\tau\to y
    &&\text{strongly in }\ L^2(I;H^{1}(\varOmega;\R^d))\,.  \label{conv4}
%L^r(I;H^{1}H^1(\varOmega;\R^d))\ \text{for all $1\leq r<\infty$} 
\end{align}\end{subequations}

We now want to pass to the limit in \eqref{momentum-weak+} as $\tau\to 0$. 
The only problematic term there is the one with $\mathscr{H}'$, because 
the other terms converge strongly due to \eqref{conv4} or are essentially linear (the dissipation term involving $\bbD$ is linear in $\nabla\DT y_\tau$, while its other factors converge strongly). We now exploit the strict monotonicity of $\mathscr{H}'$
to obtain better convergence for $\nabla^2\bary_\tau$. Consider the test functions
\[
  z_\tau:=\phi\cdot (y-\bary_\tau),\quad \text{where $\phi\in L^\infty(I;C^2(\barOmega;\R^+))$ such that 
	$\int_{\varGamma}\phi \,\d \barsigma_\tau =0$ on $I$ for all $\tau$}.
\]
Notice that $z_\tau\to 0$ strongly in $L^2(I;H^{1}(\varOmega;\R^d))$ and weakly in $L^p(I;W^{2,p}(\varOmega;\R^d))$ by \eqref{conv4} and \eqref{conv1}.
Using \eqref{momentum-weak+}, we get that
\begin{align*}%\nonumber
%&\epsilon \limsup_{\tau\to 0} \int_Q \phi |\nabla^2y-\nabla^2\bary_\tau|^p\,\d x\d t \\
	&\limsup_{\tau\to 0} \int_Q \phi[\mathscr{H}'(\nabla^2 y)-\mathscr{H}'(\nabla^2\bary_\tau)]\Vdots [\nabla^2 (y-\bary_\tau)] \,\d x\d t \\
	& \ 
	\begin{aligned}[t] 
		= \ & \limsup_{\tau\to 0}
		\int_Q [\mathscr{H}'(\nabla^2 y)-\mathscr{H}'(\nabla^2\bary_\tau)]\Vdots \nabla^2 z_\tau \,\d x\d t \\
		= \ & \limsup_{\tau\to 0} 
		\begin{aligned}[t]
			\int_Q \Big(-\bbD(\nabla\underline y_\tau)([\nabla\DT
			y_\tau]^\top\nabla\underline{y}_\tau{+}[\nabla\underline{y}_\tau]^\top\nabla\DT y_\tau)
			{:}([\nabla\underline{y}_\tau]^\top\nabla z_\tau+[\nabla z_\tau]^\top\nabla \underline{y}_	\tau) & \\
			-\varphi'(\nabla\bary_\tau){:}\nabla z_\tau +\int_Q\barf_\tau{\cdot}z_\tau & \Big) \,\d x\d t=0\,.
		\end{aligned}\\
		%= \ & 0.
	\end{aligned}
\end{align*}
Since $\phi\in L^\infty(I;C^2(\barOmega;\R^+))$ was arbitrary apart from the
requirement that $\int_{\varGamma}\phi \,\d \barsigma_\tau =0$ on $I$ for all
(small enough) $\tau$,
the strict monotonicity \eqref{ass-H2} of $\mathscr{H}'$ thus implies that 
\begin{align}\label{ytau-strong-conv}
  \nabla^2 y_\tau\to \nabla^2 y\ \ \
  \text{in $L^p(I;L^p(\varOmega{\setminus}U_t;\R^d))$, \ \ i.e., } 
	\int_I\int_{\varOmega\setminus U_t} \! |\nabla^2 y_\tau{-}\nabla^2 y|^p \,\d x \d t\to 0,
\end{align}
for any measurable set $U \subset I\times\R^d$, $U=\bigcup_{t\in I} \{t\}\times U_t$, such that 
\begin{align}\label{U-properties}
&\begin{aligned}
	\text{the support of $\barsigma_\tau(t)$ is contained in the interior of $U_t$}&\\
	\text{ for a.e.~$t\in I$ and (small enough) $\tau>0$}&.
\end{aligned}
\end{align}
In particular, $y_\tau\to y$ strongly in $L^p(I;W_{\rm loc}^{2,p}(\varOmega;\R^d))$ because $U:=I\times V$ is admissible for any closed neighborhood $V$ of $\varGamma$ in $\R^d$.

In view of \eqref{U-properties} and \eqref{ytau-strong-conv}, it is clear that for a.e.~$t\in I$, the limit function $y$ 
inherits its trace on the boundary as a strong limit of the traces of $y_\tau$, except on the part 
$\varGamma^*_t$ of $\varGamma$ that is always excluded by $U_t$, that is,
\begin{align}\label{limit-self-contact}
	\varGamma^*_t:=\left\{x\in \GNeu \,\left| \,	
	\begin{aligned}[c]
		\text{$x=\lim_{n\to \infty} x_{\tau(n)}$ for a suitable subsequence $(\tau(n))$ of $(\tau)$}&\\[-.5em]
		\text{and points $x_{\tau(n)} \in \operatorname{supp} \barsigma_{\tau(n)}(t)\subset \GNeu$}&
	\end{aligned}
	\right. \right\}\subset \scS_t.
\end{align}
Here, we indeed have that $\varGamma^*_t\subset \scS_t$ (the self-contact set of $y$ at time $t$, cf.~Definition~\ref{def}), due to
\eqref{PalHea-traction}, the definition of $\varGamma^*_t$ and the fact that $\bary_\tau(t,\cdot)\to y(t,\cdot)$ strongly in $C(\barOmega;\R^d)$, 
the latter by \eqref{conv3} and compact embedding.
Moreover,  
as a consequence of \eqref{PalHea-traction}, \eqref{conv3} and \eqref{U-properties},
\begin{align}\label{limtraction1}
	\text{$\varGamma^*_t$ contains the support of $\mathfraks_1(t,\cdot)\in W^{2-1/p,p}(\GNeu;\R^d)^*$.}
\end{align}
As a consequence of \eqref{ytau-strong-conv} and \eqref{ass-H1},
\begin{align}\label{ytau-strong-conv+}
  \mathscr{H}'(\nabla^2\bary_\tau) \to \mathscr{H}'(\nabla^2 y) \ \text{strongly in 
    $L^{p'}(I; L^{p'}(\varOmega{\setminus}U_t;\R^d))$.}
\end{align}
On the other hand, by \eqref{est-y}, passing to a subsequence if necessary,
there exists \\ $\mathfrakh\in L_{\rm w}^\infty(I; W^{2,p}(\varOmega;\R^d)^*)$ such that
\begin{align}\label{DHy-limit}
	\mathscr{H}'(\nabla^2\bary_\tau) \to \mathfrakh \ \text{weakly* in $L_{\rm w}^\infty(I; W^{2,p}(\varOmega;\R^d)^*)=L^1(I; W^{2,p}(\varOmega;\R^d))^*$}.
\end{align}
Interpreting $\mathscr{H}'(\nabla^2 y)$ as a distribution in $L^1(I; W^{2,p}(\varOmega;\R^d))^*$, i.e.,
\[
	\langle \mathscr{H}'(\nabla^2 y), z \rangle:=\int_Q \mathscr{H}'(\nabla^2 y) \Vdots \nabla^2 z \,\d x \d t
	\quad\text{for $z\in L^1(I; W^{2,p}(\varOmega;\R^d))$},
\]
we see that 
\begin{align}\label{DHy-limit-properties}
	\langle \mathfrakh, z \rangle=\langle \mathscr{H}'(\nabla^2 y), z \rangle
	\quad\text{$\forall\,z\in L^1(I; W^{2,p}(\varOmega;\R^d))$ with $z=0$ on $\varSigma^*:=\bigcup_{t\in I}\big(\{t\}\times \varGamma^*_t\big)$}.
\end{align}
In particular,
\begin{align}\label{limtraction2}
\begin{aligned}
	&\mathfraks_2:=\mathfrakh-\mathscr{H}'(\nabla^2 y)\in L^1(I; W^{2,p}(\varOmega;\R^d))^*
	\quad \text{is supported in }\ \varSigma^*, \ \text{and}\\
	&\mathfraks_2 \in L^1(I; W^{2-\frac{1}{p},p}(\varGamma_N;\R^d))^*=L_{\rm w}^\infty(I;W^{2-\frac{1}{p},p}(\varGamma_N;\R^d)^*).
\end{aligned}
\end{align}
The latter holds because due to \eqref{DHy-limit-properties}, $\langle \mathfraks_2,z \rangle$ actually only depends on the traces of $z(t,\cdot)$ on $\varGamma^*_t\subset \varGamma_N$, $t\in I$.
Altogether, we can now pass to the limit in \eqref{momentum-weak+}, 
using \eqref{DHy-limit}, \eqref{conv3} and
\eqref{conv4}. This yields the limit equation
\eqref{momentum-weak-def} with 
\[
\mathfraks:=\mathfraks_1+\mathfraks_2\in
L^2(I; W^{2-\frac{1}{p},p}(\varGamma_N;\R^d)^*).
\]
Moreover, by \eqref{limit-self-contact}, \eqref{limtraction1}, and
\eqref{limtraction2}, the total contact reaction force $\mathfraks$ satisfies
\eqref{PalHea-traction-limit}, and we conclude that $(y,\mathfraks)$ is a weak
solution to the initial-boundary-value problem
\eqref{momentum-eq}--\eqref{BC}--\eqref{IC} in
the sense of Definition~\ref{def}.
%\end{proof}
\hfill$\Box$

 \begin{remark}[Open problem: actual reaction force]\upshape
\COL{%
   In the static situations, the contact reaction force $\mathfraks$
   is a measure as shown in \cite{PalHea17ISCS}, 
and it has a natural pullback to the reference configuration given by
   $\vec n\sigma=(\nabla y|_{\GNeu}^{})^{\top}\mathfraks$.
   By contrast, in the evolution case, we are loosing this property, cf.\
   the estimate \eqref{est-traction}, because $\nabla y|_{\SNeu}^{}$ is not regular enough to justify multiplication
   with the distribution $\mathfraks$.
   Additional information about the reaction force in the static case contained
   in the ``complementary slackness principle'' of \cite{PalHea17ISCS} is also
   lost in the limit, as least if 
   $\mathfraks_2$ does not vanish.
	}
 \end{remark}

 \begin{remark}[Open problem: dynamical problems]\upshape
In many applications, inertia cannot be neglected. Yet, 
when inertial forces of the form $\varrho\DDT y$, with $\varrho>0$ a mass density in
the reference configuration, would be involved in \eqref{momentum-eq},
serious difficulties would occur in \eqref{est-traction-}
where now only the sum
$\int_{\SNeu}(\nabla\bary_\tau)^{{-}\top}\vec{n}{\cdot}z\,\d \barsigma_\tau
+\int_Q\varrho\DDT y_\tau{\cdot}z\,\d x\d t$ could be estimated. 
As a result, in the limit problem, one could not distinguish between inertial
forces and reaction forces arising from the possible self-contact. Thus, one
would have to devise a very weak solution concept.
    \end{remark}

\section*{Acknowledgements}
%         ~~~~~~~~~~~~~~~~
This research has been supported by the Czech Science Foundation
through the grants 17-04301S %(TR,
(in particular concerning %the focus on
dissipative evolutionary systems),
%18-03834S (in particular concerning the modelling of shape-memory alloys)
19-29646L
%(SK and TR,
(especially pertaining %the focus on the 
large strains in materials science)
and 19-04956S
%(TR,
(in particular concerning %the focus on
nonlinear behavior of structures).
Also the institutional support RVO:61388998 is acknowledged.

\bibliographystyle{plain}
%\bibliography{Bibliography}
\bibliography{sk-tr-viscoelasto-large}

\end{document}